\documentclass{amsart}
\usepackage[latin1]{inputenc}

\usepackage{subfigure}
\usepackage{pifont, xcolor, tikz}
\usepackage{amsmath}
\usepackage{fullpage}

\newtheorem{theorem}{Theorem}

\newtheorem{prop}[theorem]{Proposition}
\newtheorem{lemma}[theorem]{Lemma}

\newtheorem{example}[theorem]{Example}

\newcommand{\mS}{\ensuremath{\mathcal{S}}}

\newcommand{\mG}{\ensuremath{\mathcal{G}}}
\newcommand{\mV}{\ensuremath{\mathcal{V}}}
\newcommand{\bz}{\ensuremath{\mathbf{z}}}

\newcommand{\bK}{\ensuremath{\mathbf{K}}}
\newcommand{\bN}{\ensuremath{\mathbf{N}}}
\newcommand{\ox}{\ensuremath{\overline{x}}}

\newcommand{\oy}{\ensuremath{\overline{y}}}

\newcommand{\sgn}{\operatorname{sgn}}

\newcommand{\bp}{\mbox{\boldmath$\rho$}}
\newcommand{\mD}{\ensuremath{\mathcal{D}}}

\newcommand{\bone}{\ensuremath{\mathbf{1}}}

\newcommand{\diag}[1]{
  \begin{tikzpicture}[scale=0.2]\makediagb{#1}\end{tikzpicture}
}

\makeatletter
\def\testbb#1{\testbb@i#1,,\@nil}%
\def\testbb@i#1,#2,#3\@nil{%
  \draw (O) --++(#1);
  \ifx\relax#2\relax\else\testbb@i#2,#3\@nil\fi}
\makeatother   

\newcommand{\makediagb}[1]{
    \coordinate (O) at (0,0); \coordinate (N) at (0,1);
    \coordinate (NE) at (1,1); \coordinate (E) at (1,0);
    \coordinate (SE) at (1,-1); \coordinate (S) at (0,-1);
    \coordinate (SW) at (-1,-1);\coordinate (W) at (-1,0);
    \coordinate (NW) at (-1,1); \coordinate (B1) at (1.2,1.2);
    \coordinate (B2) at (-1.2,-1.2);
    
    \draw (B1) --++(0,-2.4); \draw (B1) --++ (-2.4,0);
    \draw (B2) --++(0,2.4);  \draw (B2) --++ (2.4,0);   
    \testbb{#1}
} 

\author{Stephen Melczer}
\author{Mark C. Wilson}

\title{Asymptotics of lattice walks via analytic combinatorics in several variables} 
\address[S. Melczer]{U. Lyon, CNRS, ENS de Lyon, Inria, UCBL, Laboratoire LIP \& Cheriton School of Computer Science, University of Waterloo, Waterloo ON Canada}
\email{smelczer@uwaterloo.ca}
\address[M. C. Wilson]{Department of Computer Science, University of Auckland, Private Bag 92019 Auckland, New Zealand}
\email{mcw@cs.auckland.ac.nz}

\title{Asymptotics of lattice walks via analytic combinatorics in several variables} 

\begin{document}

\begin{abstract}
\textbf{Abstract.} 
We consider the enumeration of walks on the two-dimensional non-negative integer lattice with steps defined by a finite set $\mS \subseteq \{\pm1,0\}^2$. Up to isomorphism there are 79 unique two-dimensional models to consider, and previous work in this area has used the \emph{kernel method}, along with a rigorous computer algebra approach, to show that 23 of the 79 models admit D-finite generating functions.  In 2009, Bostan and Kauers used Pad{\'e}-Hermite approximants to guess differential equations which these 23 generating functions satisfy, in the process guessing asymptotics of their coefficient sequences.  In this article we provide, for the first time, a complete rigorous verification of these guesses.  Our technique is to use the kernel method to express 19 of the 23 generating functions as diagonals of tri-variate rational functions and apply the methods of analytic combinatorics in several variables (the remaining 4 models have algebraic generating functions and can thus be handled by univariate techniques).  This approach also shows the link between combinatorial properties of the models and features of its asymptotics such as asymptotic and polynomial growth factors.  In addition, we give expressions for the number of walks returning to the $x$-axis, the $y$-axis, and the origin, proving recently conjectured asymptotics of Bostan, Chyzak, van Hoeij, Kauers, and Pech.
\\

\textbf{R\'esum\'e.} 
Nous prouvons les comportements asymptotiques conjectur{\'e}s par Bostan et Kauers pour les 23 mod{\'e}les de marches dirig{\'e}es restreintes au quadrant positif des r{\'e}seaux bi-dimensionnels. Notre technique combine la m{\'e}thode du noyau avec les outils de la combinatoire analytique en plusieurs variables. Nous prouvons aussi les comportements asymptotiques des nombres de marches terminant sur l'un des axes, ou {\`a} l'origine.
\end{abstract}

\maketitle

%-------------------------------------------------
\section{Introduction}
\label{sec:introduction}
%-------------------------------------------------

Recently, the study of two-dimensional lattice walks restricted to the non-negative quadrant has been an active topic of interest in several sub-areas of combinatorics (see, for instance,~\cite{ Bous02,BoPe03,Bous05, Mish09, MiRe09, KaKoZe09, BoKa09, BoKa10,BoMi10, Rasc12,FaRa12, KuRa12,BoRaSa14,BoKuRa14, Bo15}), with applications in branches of applied mathematics including queuing theory and the study of linear polymers.  The seminal work of Mishna and Bousquet-M{\'e}lou~\cite{ BoMi10} gave a uniform approach to several enumerative questions, including the nature of a model's GF\footnote{We abbreviate `generating function' as GF throughout.} (algebraic, D-finite, etc.) and the determination of exact or asymptotic counting formulas.  In particular, they used the \emph{kernel method} to prove that the GFs corresponding to 22 of the 79 non-equivalent two-dimensional models are D-finite. They conjectured that one additional model was D-finite---proved later by several authors~\cite{BoKa10, BoKuRa14, Bo15}---and that the rest were not.  In 2009, Bostan and Kauers~\cite{ BoKa09} used computer algebra approaches to guess differential equations satisfied by the GFs of these 23 models, which were then exploited to guess dominant asymptotics for the number of walks of a given length; see Table~\ref{tab:guessed}.  

\begin{table}[ht]
\centering
\begin{tabular}{ | c | c | c @{ \hspace{0.01in} }@{\vrule width 1.2pt }@{ \hspace{0.01in} }  c | c | c @{ \hspace{0.01in} }@{\vrule width 1.2pt }@{ \hspace{0.01in} }  c | c | c |  }
  \hline
   \# & $\mS$ & Asymptotics & \# & $\mS$ & Asymptotics & \# & $\mS$ & Asymptotics \\ \hline
  &&&&&&&& \\[-5pt] 
  \small{1} & \diag{N,S,E,W}  & $\frac4\pi \cdot \frac{4^n}n$ &
  \small{9} & \diag{NE,NW,S}  & $\frac{\sqrt{3}}{2\sqrt{\pi}} \cdot \frac{3^n}{\sqrt{n}}$ &
  \small{17} & \diag{N,SE,SW} & $\frac{4 \cdot A_n}{\pi} \cdot \frac{(2\sqrt{2})^n}{n^2}$  \\
  %%%%%%%%%%%%%
  \small{2} & \diag{NE,SE,NW,SW}  & $\frac2\pi \cdot \frac{4^n}n$ &
  \small{10} & \diag{N,NW,NE,S}  & $\frac4{3\sqrt{\pi}} \cdot \frac{4^n}{\sqrt{n}}$ &
  \small{18} & \diag{N,S,SE,SW} & $\frac{3\sqrt{3}\cdot B_n}{\pi} \cdot \frac{(2\sqrt{3})^n}{n^2}$ \\
  %%%%%%%%%%%%%
  \small{3} & \diag{N,S,NE,SE,NW,SW} & $\frac{\sqrt{6}}\pi \cdot \frac{6^n}n$  &
  \small{11} & \diag{NE,NW,E,W,S}  & $\frac{\sqrt{5}}{2\sqrt{2\pi}} \cdot \frac{5^n}{\sqrt{n}}$ &
  \small{19} & \diag{N,E,W,SE,SW}  & $\frac{\sqrt{8}(1+\sqrt{2})^{7/2}}{\pi} \cdot \frac{(2+2\sqrt{2})^n}{n^2}$ \\
  %%%%%%%%%%%%%
  \small{4} & \diag{N,S,E,W,NW,SW,SE,NE}  & $\frac{8}{3\pi} \cdot \frac{8^n}n$ &
  \small{12} & \diag{N,NE,NW,SE,SW}  & $\frac{\sqrt{5}}{3\sqrt{2\pi}} \cdot \frac{5^n}{\sqrt{n}}$ &
  \small{20} & \diag{NE,NW,SE,SW,S}  & $\frac{6C_n}{\pi} \cdot \frac{(2\sqrt{6})^n}{n^2}$ \\
  %%%%%%%%%%%%%
  \small{5} & \diag{NE,W,S}   & $\frac{2\sqrt{2}}{\Gamma(1/4)} \cdot \frac{3^n}{n^{3/4}}$ &
  \small{13} & \diag{N,NW,NE,E,W,S} & $\frac{2\sqrt{3}}{3\sqrt{\pi}} \cdot \frac{6^n}{\sqrt{n}}$ &
  \small{21} & \diag{N,E,W,S,SW,SE}  & $\frac{\sqrt{3}(1+\sqrt{3})^{7/2}}{2\pi} \cdot \frac{(2+2\sqrt{3})^n}{n^2}$  \\
  %%%%%%%%%%%%%
  \small{6} & \diag{N,E,SW}   & $\frac{3\sqrt{3}}{\sqrt2\Gamma(1/4)} \cdot \frac{3^n}{n^{3/4}}$ &
  \small{14} & \diag{N,E,W,NE,NW,SE,SW} & $\frac{\sqrt{7}}{3\sqrt{3\pi}} \cdot \frac{7^n}{\sqrt{n}}$ & 
  \small{22} & \diag{NE,NW,E,W,SE,SW,S}  & ${\scriptstyle \frac{\sqrt{6(379+156\sqrt{6})(1+\sqrt{6})^7}}{5\sqrt{95}\pi} \cdot \frac{(2+2\sqrt{6})^n}{n^2}}$ \\
  %%%%%%%%%%%%%
  \small{7} & \diag{N,NE,E,S,SW,W}   & $\frac{\sqrt{6\sqrt{3}}}{\Gamma(1/4)} \cdot \frac{6^n}{n^{3/4}}$ &
  \small{15} & \diag{N,W,SE}   & $\frac{3\sqrt{3}}{2\sqrt{\pi}} \cdot \frac{3^n}{n^{3/2}}$ &
  \small{23} & \diag{E,SE,W,NW}  & $\frac{8}{\pi} \cdot \frac{4^n}{n^2}$ \\
  %%%%%%%%%%%%%
  \small{8} & \diag{NE,E,SW,W}  & $\frac{4\sqrt3}{3\Gamma(1/3)} \cdot \frac{4^n}{n^{2/3}}$ &
  \small{16} & \diag{NW,SE,N,S,E,W}   & $\frac{3\sqrt{3}}{2\sqrt{\pi}} \cdot \frac{6^n}{n^{3/2}}$  &&& \\
\hline
\end{tabular}
\vspace{0.05in}

\caption{Asymptotics for the 23 D-finite models.} \label{tab:guessed}
\vspace{-0.4in}

\[ {\scriptstyle A_n = 4(1-(-1)^n)+3\sqrt{2}(1+(-1)^n), \quad B_n = \sqrt{3}(1-(-1)^n)+2(1+(-1)^n), \quad C_n = 12/\sqrt{5}(1-(-1)^n)+\sqrt{30}(1+(-1)^n)} \]
\end{table} 

The main difficulty one encounters in trying to determine asymptotics through annihilating equations of D-finite GFs is the \emph{connection problem}; that is, even if one is able to rigorously derive an annihilating linear differential equation and compute asymptotics for a basis of solutions, it can be surprisingly hard (possibly incomputable in general) to rigorously express the GF in question as a linear combination of these basis elements.  Ongoing work of Bostan, Chyzak, van Hoeij, Kauers, and Pech~\cite{BCHKP15+} attempts to get around this problem by using creative telescoping techniques combined with the kernel method to represent the walk GFs explicitly in terms of hypergeometric functions.  Although such a representation should, in principle, allow one to rigorously determine asymptotics, in practice this depends on computing integrals of hypergeometric functions which those authors have only been able to numerically approximate\footnote{For some models, such integrals need to be rigorously determined to show not only the asymptotic constant of growth but even its exponential growth.}.

\subsection{Our Contribution}  In this work, we combine the expressions resulting from the kernel method with the techniques of analytic combinatorics in several variables to rigorously determine the asymptotics for 19 of the 23 D-finite models, verifying the computational guesses of \cite{BoKa09}\footnote{Three of the models have a periodic constant which their original table failed to take into account, but their guesses are otherwise correct.}.  The final 4 models admit algebraic GFs whose minimal polynomials are explicitly known, meaning their asymptotics can be determined rigorously through univariate means\footnote{While these algebraic equations can be used to derive diagonal expressions for the 4 algebraic GFs, the resulting representations have more pathological characteristics (such as degenerate critical points) than those arising `naturally' out of the kernel method.  In any case, univariate methods for algebraic GFs are well established and easily give asymptotics for these models.} \cite{ Mish09,Bous10, BoKa10}.  Thus, to our knowledge, this work gives the first complete proof of asymptotics for the 23 D-finite models.

The analysis breaks down into four groups (up to exchanging $x$ and $y$ coordinates): 4 models whose step sets are symmetric over every axis (the \emph{highly symmetric models}, these were handled in general dimension by Melczer and Mishna~\cite{ MeMi15}), 6 models whose step sets are symmetric over the $y$ axis and have positive vector sum in their second coordinate (the \emph{positive drift} models), 6 models whose step sets are symmetric over the $y$ axis and have negative vector sum in their second coordinate (the \emph{negative drift} models), and 3 sporadic cases.  One reason for the recent interest in lattice path models in the quarter plane is the large variety of asymptotic and analytic behaviour their GFs exhibit, and this is evident in our work as well.  The GF of a model is represented as the diagonal of a multivariate rational function, and the range of asymptotic behaviour apparent in Table~\ref{tab:guessed} reflects differences in the geometry of the set of singularities of these rational functions.

We begin in Section~\ref{sec:KM} by giving an overview of lattice path enumeration and the kernel method, and show how it can be used to derive expressions for lattice path GFs which are amenable to the techniques of analytic combinatorics in several variables.  Section~\ref{sec:ACSV} then details the general methods of analytic combinatorics in several variables, and outlines how the asymptotic analysis will proceed.  This is followed by Section~\ref{sec:Asm}, where we derive asymptotics for the 19 models represented by multivariate diagonals through the kernel method.  In Section~\ref{sec:Boundary} we examine the asymptotics for walks which return to one or both of the boundaries $x=0$ and $y=0$---proving recently guessed asymptotics by Bostan et al.~~\cite{BCHKP15+}---and conclude in Section~\ref{sec:Conc} with directions for future research.  For more historical background and details on the calculations, the reader is referred to an upcoming full journal article building upon this extended abstract.

\section{The Kernel Method for Quadrant Walks}
\label{sec:KM}

The kernel method is a widely used strategy for manipulating functional equations, often those arising in the context of enumerating lattice paths in restricted regions.  An often cited early example appears in the work of Knuth~\cite{Knu68}, and more modern accounts include Bousquet-M{\'e}lou and Petkov\v{s}ek~\cite{ BoPe00} and Banderier et al.~\cite{BaBoAll02}, among many others.

To a given set of steps $\mS = \{\pm1,0\}^2\setminus\{(0,0)\}$ we associate the multivariate GF
\[ C(x,y,t) := \sum_{i,j,n \geq 0} c_{i,j,n}x^iy^jt^n, \]
where $c_{i,j,n}$ denotes the number of walks on the steps $\mS$ of length $n$, beginning at the origin, ending at the lattice point $(i,j)$, and never leaving the first quadrant (including coordinate axes).  We further define the \emph{characteristic polynomial} of a model to be the Laurent polynomial $S(x,y) = \sum_{(i,j)\in\mS}x^iy^j$. Decomposing a walk of length $n>0$ ending at the point $(i,j)$ into a walk of length $n-1$ followed by a step in $\mS$ gives a recurrence for $c_{i,j,n}$ (with special cases when $i=0$ or $j=0$ to account for the restriction to the quarter plane).  Translating this recurrence into GF equations allows one to show that the GF $C(x,y,t)$ satisfies the functional equation
\[ C(x,y,t) = 1 + tS(x,y)C(x,y,t) - t\ox([\ox]S(x,y)) C(0,y,t) - t\oy([\oy]S(x,y)) C(x,0,t) + \epsilon tC(0,0,t),\]
where $\ox$ is shorthand for $1/x$ and $\epsilon=1$ if $(-1,-1)\in\mS$ and 0 otherwise.  Letting $B_{-1}(y)=[\ox]S(x,y)$ and $A_{-1}(x) = [\oy]S(x,y)$, this is typically written in the form
\begin{equation} xy(1-tS(x,y))C(x,y,t) = xy - tyB_{-1}(y)C(0,y,t) - txA_{-1}(x)C(x,0,t) + \epsilon xyt C(0,0,t), \label{eq:funC} \end{equation}
and we note that $C(0,y,t),C(x,0,t),$ and $C(0,0,t)$ represent the GFs for walks ending on the $x$-axis, on the $y$-axis, and at the origin, respectively.
\\

In order to manipulate this functional equation into a more usable form, Bousquet-M{\'e}lou~\cite{Bous05} (see also \cite{BoMi10}) 
%borrowed from \cite{FaIaMa99} the idea of introducing 
used a group of bi-rational transformations of the plane which fix $S(x,y)$.  As $\mS \subset \{\pm1,0\}^2 \setminus \{(0,0)\}$, we can in fact write
\[ S(x,y) = \oy A_{-1}(x) + A_0(x) + y A_1(x) = \ox B_{-1}(y) + B_0(y) + x B_1(y),  \]
for Laurent polynomials $A_i,B_i$.  If $A_{-1}(x)=0$ or $B_{-1}(y)=0$ then $\mS$ has no step moving towards (at least) one of its boundaries.  In other words, the model defined by $\mS$ is actually a lattice path model with a restriction to a halfplane (or having no restriction).  Banderier and Flajolet~\cite{ BaFl02} have shown that such models always admit algebraic GFs, and gave effective means of calculating their asymptotics.  Thus, we may assume neither $A_{-1}(x)$ nor $B_{-1}(y)$ is 0.  With this in mind, we define the bi-rational transformations
\[ \Phi:(x,y) \rightarrow \left(x, \oy \frac{A_{-1}(x)}{A_1(x)}\right) \qquad \Psi:(x,y) \rightarrow \left(\ox \frac{B_{-1}(y)}{B_1(y)}, y \right), \]
and let $\mG$ be the (possibly infinite) group of bi-rational transformations these involutions generate. We can view an element $g \in \mG$ as acting on a Laurent polynomial $f(x,y)$ by setting $\sigma(f(x,y)) := f(\sigma(x,y))$. 
\\

Define the non-negative series extraction operator $[x^{\geq}y^{\geq}]:\mathbb{Q}[x,\ox,y,\oy][[t]]\rightarrow\mathbb{Q}[x,y][[t]]$ by 
\[ [x^{\geq}y^{\geq}] \sum_{n \geq 0}\left( \sum_{i,j \in \mathbb{Z}}r_{i,j,n}x^iy^j \right)t^n := \sum_{n,i,j \geq 0} r_{i,j,n}x^iy^jt^n. \]
The main results of Bousquet-M{\'e}lou and Mishna, and Bostan and Kauers, combine to yield the following.

\begin{theorem}[\cite{ BoMi10, BoKa10}]
\label{thm:BoMi}
Up to isomorphism, there are 79 distinct two-dimensional lattice path models with short steps restricted to the non-negative quadrant (which are not equivalent to halfplane models).  Of these 79, precisely the 23 models listed in Table~\ref{tab:guessed} give rise to a finite group $\mG$.  For each model in Table~\ref{tab:guessed} except models 5--8, the multivariate GF $C(x,y,t)$ can be expressed as
\begin{equation} C(x,y,t) =  [x^{\geq}y^{\geq}] \frac{\sum_{g \in \mG} \sgn(g) g(xy)}{xy(1-tS(x,y))}, \label{eq:OS} \end{equation}
where $\sgn(g)$ is the minimal length of $g$ in terms of the generators $\Phi$ and $\Psi$. 
For models 5--8, $C(x,y,t)$ is algebraic and the dominant asymptotics of $C(1,1,t)$ match Table~\ref{tab:guessed}.
\end{theorem}

The key idea is that $\Phi$ and $\Psi$ each fix $S(x,y)$ along with one term on the right hand side of Equation~\eqref{eq:funC}.  Thus the sign-weighted sum over the group eliminates all unknown series on the right hand side.  This introduces new terms $C(g(x,y),t)$ on the left hand side for each $g \in \mG$, but one can show that for each case, except models 5--8, taking the non-negative series extraction leaves only $C(x,y,t)$.  The algebraic models 5, 6, and 7 require a more delicate analysis \cite{BoMi10}, and model 8 a different approach \cite{BoKa10}. 

The 56 models with infinite group are strongly suggested to have non-D-finite GFs $C(1,1,t)$ counting the total number of walks, meaning that determining their asymptotics is likely to be difficult.  Dominant asymptotics are only known for 5 of the 51 models (see Melczer and Mishna~\cite{ MeMi14a}) although Bostan et al.~\cite{ BoRaSa14} have determined asymptotics up to a constant multiple for the number of walks returning to the origin in the other 51 cases. Recent work of Duraj~\cite{Du14} implies that when the vector sum of a model's step set contains two negative coordinates the results of~\cite{BoRaSa14} determine asymptotics up to a constant multiple for the number of walks in these models ending anywhere.  Fayolle and Raschel~\cite{FaRa12} outline a method which one can use to find the exponential growth of all of these 51 models, although they cannot generally recover sub-exponential behaviour.

\subsection{Diagonal expressions}
In order to prove the dominant asymptotics listed in Table~\ref{tab:guessed}, we will need to convert Equation~\eqref{eq:OS} into a more computational form.  For that we define the diagonal operator $\Delta:\mathbb{Q}[x,\ox,y,\oy][[t]]\rightarrow\mathbb{Q}[[t]]$ by 
\[ \Delta  \left(\sum_{n \geq 0}\left( \sum_{i,j \in \mathbb{Z}}r_{i,j,n}x^iy^j \right)t^n \right) := \sum_{n \geq 0} r_{n,n,n}t^n. \]
The result which allows us a compact representation of the GF $C(1,1,t)$ for the number of walks ending anywhere---along with the GFs $C(1,0,t),C(0,1,t)$ and $C(0,0,t)$ for walks ending on one or both of the coordinate axes---is the following.

\begin{lemma} \label{lem:diagsubs}
Let $P(x,y,t) \in \mathbb{Q}[x,\ox,y,\oy][[t]]$ and $Q(x,y,t) = [x^\geq y^\geq]P(x,y,t)$.  Then 
\begin{align*} \label{eq:postodiag}
Q(1,1) &= \Delta R(x,y,t) & 
Q(0,1) &= \Delta\left((1-x)\cdot R(x,y,t)\right) \\
Q(1,0) &= \Delta\left((1-y)\cdot R(x,y,t)\right) &
Q(0,0) &= \Delta\left((1-x)(1-y)\cdot R(x,y,t)\right),
\end{align*}
where 
\[ R(x,y,t) = \frac{P\left(\ox,\oy,xyt\right)}{(1-x)(1-y)}. \]
\end{lemma} 

The proof follows from the definition of the diagonal after writing out the geometric series and expansion of $P$ on the right hand side.  Combining Lemma~\ref{lem:diagsubs} with Theorem~\ref{thm:BoMi} gives us a compact representation for our GFs which will allow for the asymptotic analysis.

\begin{theorem} \label{thm:diag}
Let $\mS$ be a step set corresponding to a model in Table~\ref{tab:guessed} other than the algebraic models 5--8.  Then for $a,b \in \{0,1\}$,
\begin{equation} C(a,b,t) = \Delta\left( \frac{xyO(\ox,\oy)}{1-txyS(\ox,\oy)} \cdot (1-x)^{-a}(1-y)^{-b} \right), \label{eq:diag} \end{equation}
where $O(x,y) = \sum_{g \in \mG}\sgn(g)g(xy)$.
\end{theorem}

%%%%%%%%%%%%%%%%
% ACSV Section
%%%%%%%%%%%%%%%%

\section{Analytic Combinatorics in Several Variables}
\label{sec:ACSV}

A general reference for this section is the text of Pemantle and Wilson~\cite{PeWi13} which contains precise statements and proofs for the results sketched here. Let $F(x,y,t)  = G(x,y,t) / H(x,y,t) \in \mathbb{Q}(x,y,t)$ be a rational function analytic at the origin.  Similar to the univariate case, the starting point here is to use the (multivariate) Cauchy integral formula to write
\begin{equation} a_{n,n,n} = \frac{1}{(2\pi i)^3}\int_T F(x,y,t)\frac{dx}{x^{n+1}}\frac{dy}{y^{n+1}}\frac{dt}{t^{n+1}}, \label{eq:CI} \end{equation}
where $T$ is a sufficiently small torus around the origin. The set of zeros $\mV = \mathbb{V}(H) \cap (\mathbb{C}^*)^3$ describes the set of singularities of $F(x,y,t)$ off the coordinate axes, and is known as the \emph{singular variety}.  If $\mD$ is the open domain of convergence of $F$ at the origin, we call a singular point $(x,y,t) \in \mV \cap \partial \mD$ on the boundary of the domain of convergence a \emph{minimal point}.  Minimal points are, in a sense, the most straightforward generalization of dominant singularities in the univariate case.  Unfortunately, however, the extreme pathology of singularities possible in the multivariate case means it is possible that no minimal points will contribute to the dominant asymptotics.  

Another notion is needed to fill this gap.  The function $h : \mV \rightarrow \mathbb{R}$ defined by $h(x,y,t) := -\log|xyt|$ is known as the \emph{height function} associated to $\mV$ and captures the part of the integrand in Equation~\eqref{eq:CI} whose modulus increases with $n$.  Assume first that $\mV$ defines a smooth manifold.  Then the work of Pemantle and Wilson implies that the critical points of the function $h$ as a map between smooth manifolds are those around which local behaviour of the GF will determine asymptotics of the coefficient sequence, they are called the \emph{critical points} of $\mV$.  Basic results from complex analysis imply that we can deform the cycle of integration $T$ in Equation~\eqref{eq:CI} without changing the value of the integral as long as we do not cross the singular variety $\mV$ or the coordinate axes.  If there are critical points on the boundary of the domain of convergence---i.e., critical minimal points---then one can deform $T$ to be arbitrarily close without changing the integral representation of $a_{n,n,n}$.  This allows one to show that the local behaviour of the function at these points will determine asymptotics (under the assumption that $\mV$ is smooth).

In the general case, one begins by computing a \emph{Whitney stratification} of the singular variety, which is a decomposition of $\mV$ into a disjoint collection of smooth manifolds called \emph{strata} with some additional properties (see Pemantle and Wilson~\cite{PeWi13}, Definition 5.4.1 and the following discussion). When $F$ is rational, a stratification can be computed algorithmically: in the smooth case the stratification simply consists of the manifold itself and in general each stratum can be effectively represented as the intersection of algebraic hypersurfaces $\mathbb{V}(H_1)\cap\cdots\cap\mathbb{V}(H_k)$ minus some varieties of lower dimension.  A point $(x,y,t)$ in a stratum $B$ is called a \emph{critical point of $B$} if $\nabla h|_B(x,y,t) = 0$, and the set of \emph{critical points} of the singular variety is the set of points that are critical points for some stratum.

When $\mV$ is a smooth manifold, or is smooth except for points where it consists of smooth manifolds intersecting transversely (the \emph{multiple point} case) the set of critical points forms an algebraic set easily computed by Gr\"{o}bner basis or other elimination methods;  generically the set of critical points is finite. It is very difficult in general to determine which critical points actually contribute to the dominant asymptotics.  

For a point $(x,y,t)$ where $\mV$ is locally the transverse intersection of smooth hypersurfaces $\mathbb{V}(H_1) \cap \cdots\cap\mathbb{V}(H_k)$ we define the cone $\bK(x,y,t) \subset \mathbb{RP}^2$ to be the span of the vectors
\[ \nabla_{\log} H_i := (x \partial H_i/\partial x, y\partial H_i/\partial y, t\partial H_i/\partial t),  \quad i=1,\dots,k \]
and the cone $\bN(x,y,t) = \bK(x,y,t)^*$ to be the dual cone to $\bK(x,y,t)$. Pemantle and Wilson~\cite{PeWi08} proved that a multiple point $(x,y,t)$ is critical if and only if $(1,1,1) \in \bN(x,y,t)$.  Furthermore, they showed that when critical minimal points exist they almost always are the ones determining asymptotics for the diagonal coefficient sequence.

\begin{prop}[Pemantle and Wilson~\cite{PeWi13}, Proposition 10.3.6] \label{prop:contrib}
Suppose that for a rational function $G/H$ the singular variety $\mV$ is composed solely of smooth or multiple points and let $K$ be its set of critical points.  Let $W \subset \overline{\mD}$ be the points in the closure of the power series domain of convergence at which $h(x,y,t)$ is minimized, and assume that the function $(x,y,t)\mapsto(|x|,|y|,|t|)$ is constant on $W$. If $\bone \notin \partial \mathbf{N}(x,y,t)$ for any $(x,y,t) \in K \cap W$ and the set
\[ V = \{ (x,y,t) \in K \cap W : (1,1,1) \in \mathbf{N}(\bz) \} \]
is non-empty and finite, then $V$ is the set of \emph{contributing} points of $F$, in the sense that
\[ [x^ny^nt^n]F(x,y,t) \sim \sum_{(x,y,t) \in V} \textsl{formula}(x,y,t), \]
where $\textsl{formula}(x,y,t)$ denotes an effective function---under mild conditions---which depends on the local geometry of $\mV$ at $(x,y,t)$.
\end{prop}
A precise description of formula$(x,y,t)$ in the smooth and multiple point cases can be found in~\cite[Thm 3.2]{RaWi11}, which is the version used in our calculations through a Sage implementation by Raichev~\cite{Raic12}.

\section{Asymptotics}
\label{sec:Asm}
We now apply the results of analytic combinatorics in several variables to the expressions obtained in Theorem~\ref{thm:diag}, proving the guesses of Bostan and Kauers~\cite{ BoKa09}.

\subsection{The Highly Symmetric Models}
\[ \diag{N,S,E,W} \qquad \diag{NE,SE,NW,SW} \qquad \diag{N,S,NE,SE,NW,SW} \qquad \diag{N,S,E,W,NW,SW,SE,NE} \]
Four of the models in Table~\ref{tab:guessed} have step sets which are symmetric over every axis, and for each model one can directly calculate that the group $\mG$ is $\{(x,y) \mapsto (x^{\pm1},y^{\pm1})\}$, meaning Equation~\eqref{eq:diag} simplifies to
\[  C(1,1,t) = \Delta\left( \frac{(1+x)(1+y)}{1-txyS(x,y)} \right). \]
Let $H=1-txyS(x,y)$.  There are no solutions to $H = H_t = 0$, so $\mV$ is smooth, and the condition $(1,1,1) \in \bK(x,y,t)$ becomes $xH_x=yH_y=tH_t$.  Solving this yields the critical points: $K = \{(\pm1,\pm1,\pm|\mS|)\} \cap \mV$.  Melczer and Mishna~\cite{ MeMi15} showed---for the analogue in dimension $d$---that the points in $K$ are all minimal.  Thus, one can use Proposition~\ref{prop:contrib} along with \cite[Thm 3.2]{RaWi11} to give the asymptotics of $C(1,1,t)$ which appear in Table~\ref{tab:guessed}.

\subsection{Models With One Symmetry}
There are 12 models whose step sets have one symmetry (we assume without loss of generality over the $y$ axis).  For each of these models, the group $\mG$ is a group of order 4 generated by $(x,y)\mapsto(\ox,y)$ and $(x,y)\mapsto(x,\oy A_{-1}(x)/A_1(x))$, and Equation~\eqref{eq:diag} simplifies to
\[ C(1,1,t) = \Delta\left(\frac{(1+x)\left(A_1(x)-y^2A_{-1}(x)\right)}{A_{1}(x)(1-y)(1-txyS(x,\oy))} \right). \]
We note that this rational function may be singular at the origin (if $A_{-1}=\ox+x$ and $A_1=1$, for instance) but if it is not analytic it is of the form $R(x,y,t)/x$ for $R$ analytic at the origin, and we can use the identity $[t^n]\Delta(R/x) = [t^{n+1}]\Delta(ytR)$ to determine the asymptotics of the diagonal sequence by analyzing the function $ytR(x,y,t)$ which is analytic at the origin.

Let $H_1=A_1(x), H_2=1-y,$ and $H_3=1-txyS(x,\oy)$.  The singular variety $\mV$ is the union of the three smooth varieties $\mV_1 = \mathbb{V}(H_1),\mV_2 = \mathbb{V}(H_2),$ and $\mV_3 = \mathbb{V}(H_3)$, which intersect transversely.  As $H_3$ is the only factor with $t$ any critical point $(x,y,t)$ must be in $\mV_3$, so there are four possible strata which could provide critical points.  The first two cases are:

\begin{enumerate} \itemsep=0.2em
\item \textbf{(critical smooth point of $\mV_3$)} As the multivariate expansion of $1/H_3(x,y,t)$ has all non-negative coefficients, $(x,y,t)$ is a critical minimal point of $\mV_3$ only if $(|x|,|y|,|t|)$ is (see~\cite{PeWi08}).  Thus, we search for positive real solutions of the smooth critical point equations $xH_x=yH_y=tH_t$.  This simplifies to $S_x(x,\oy)=S_y(x,\oy)=0$, and due to the step set symmetry over the $x$ axis one finds the only positive real solution occurs at $(x_1,y_1,t_1) = \left(1,\sqrt{A_1(1)/A_{-1}(1)}, 1/(x_1y_1S(x_1,\oy_1)) \right).$
\item \textbf{(critical multiple  point of $\mV_2 \cap \mV_3 \setminus \mV_1$)} Any point $(x,y,t)$ on $\mV_2$ has $y=1$ and, as $\Delta_{log}(H_2) = (0,-1,0)$ at any point on $\mV_2$, we see $(1,1,1) \in \bK(x,y,t)$ if and only if $S_y(x,1)=1$.  This gives one positive real critical point in this stratum: $(x_2,y_2,t_2) = (1,1,1/|\mS|)$.
\end{enumerate}  
The remaining two strata which could contribute critical points are $\mV_1\cap\mV_3 \setminus \mV_2$ and $\mV_1\cap\mV_2\cap\mV_3$, however they do not.  If $A_1(x)=0$ and $\mS$ is not symmetric over both axes, a straightforward computation shows $S_y(x,\oy)=0$ has no solution when $y\neq0$.  
%Thus, there are two possible positive critical minimal points: a smooth point on $\mV_3$ and a multiple point on $\mV_2 \cap \mV_3$.  

As $t=1/(xyS(x,\oy))$ on $\mV_3$, any positive real point $(x,y,t) \in \mV$ can be shown to be minimal with respect to $\mV_3$.  After simplification, $H_1(x)=1,1+x^2,$ or $1+x+x^2$ so the factor $H_1$ does not affect the minimality of the two points above.  Whether or not 
$\mV_2$ affects the analysis depends on whether $A_1(1) < A_{-1}(1)$ or $A_1(1) > A_{-1}(1)$ (if they are equal we are in the highly symmetric case).  We now examine both of these cases.

\subsubsection{The Negative Drift Models}
\label{sec:posdrift}
\[ \diag{N,SE,SW} \qquad \diag{N,S,SE,SW} \qquad \diag{NE,NW,SE,SW,S} \qquad \diag{N,E,W,SE,SW} \qquad \diag{N,E,W,S,SW,SE}  \qquad \diag{NE,NW,E,W,SE,SW,S} \]

If $A_1(1) < A_{-1}(1)$---i.e., more steps move south than north---then both points we have found above are critical minimal points. One can verify that $(x_1,y_1,t_1)$ minimizes the height function---as the minimum must occur at a critical point --- so the points contributing to the dominant asymptotics are $(x_1,y_1,t_1)$ and any other points on $\mV_3$ with the same modulus.  It can be shown that the contributing points are exactly the set 
\[ \{(x,y,t) : x=\pm1,\quad y^4=\sqrt{A_1(x)/A_{-1}(x)},\quad t= 1/xyS(x,\oy),\quad |y|=|y_1|,\quad |t| = |t_1| \}, \]
and each is a smooth point of the singular variety.  Theorem 3.2 of Raichev and Wilson~\cite{RaWi11} calculates the contribution of each point and gives the asymptotics of $C(1,1,t)$ which appear in Table~\ref{tab:guessed}.

\begin{example} Consider the model defined by step set $\mS = \{ (0,1),(-1,-1),(1,-1) \} = \{N,SE,SW\}$.  Here we have 
{\small \[ [t^n]C(1,1,t) = [t^n]\Delta \left( \frac{(1+x) \left(1 - y^2 (\ox + x)\right)}{(1-y)(1-t(x + y^2+x^2y^2))} \right) = 
[t^{n+1}]\Delta \left( \frac{yt(1+x) \left(x - y^2 (x^2 + 1)\right)}{(1-y)(1-t(x + y^2+x^2y^2))} \right), \]}
and four of the eight possible points described above are contributing points: 
\[ \bp_1 = (1,1/\sqrt{2},1/2) \quad \bp_2 = (1,-1/\sqrt{2},1/2) \quad \bp_3 = (-1,i/\sqrt{2},-1/2) \quad \bp_4 = (-1,-i/\sqrt{2},-1/2).\]
Using the Sage implementation of Raichev~\cite{Raic12}, we calculate the contribution at each to be
{\small\[ 
\Psi^{(\bp_1)}_{n} = \frac{4(3\sqrt{2}+4)}{\pi} \cdot \frac{(2\sqrt{2})^n}{n^2}\qquad
 \Psi^{(\bp_2)}_{n} = \frac{4(3\sqrt{2}-4)}{\pi} \cdot \frac{(-2\sqrt{2})^n}{n^2} \qquad
 \Psi^{(\bp_3)}_{n} = \Psi^{(\bp_4)}_{n} = O\left((2\sqrt{2})^n n^{-3}\right) 
\]}
so that the number of walks of length $n$ satisfies
\begin{align*} 
s_n &= \frac{(2\sqrt{2})^n}{n^2} \cdot \frac{4}{\pi}\left(4(1-(-1)^n) + 3\sqrt{2}(1+(-1)^n)  + O(n^{-1})\right) \\
&= 
\left\{ \begin{array}{ll}    
\frac{(2\sqrt{2})^n}{n^2} \cdot \left( \frac{24\sqrt{2}}{\pi} + O(n^{-1})\right) &: n \text{ even} \\ 
\frac{(2\sqrt{2})^n}{n^2} \cdot \left( \frac{32}{\pi} + O(n^{-1})\right) &: n \text{ odd} \end{array}  \right.
\end{align*}
Note that the original table of Bostan and Kauers~\cite{ BoKa09} only had the constant for the even cases listed.
\end{example}

\subsubsection{The Positive Drift Models}
\label{sec:negdrift}
\[ \diag{NE,NW,S} \qquad \diag{N,NW,NE,S} \qquad \diag{N,NE,NW,SE,SW} \qquad \diag{NE,NW,E,W,S} \qquad   \diag{N,NW,NE,E,W,S} \qquad \diag{N,E,W,NE,NW,SE,SW} \]

When $A_1(1) > A_{-1}(1)$---i.e., more steps move north than south---then the factor $H_2=1-y$ makes $(x_1,y_1,t_1)$ fall outside the domain of convergence.  Thus, there is a single positive critical minimal point: $(1,1,1/|\mS|)$.  Now this point will determine the dominant asymptotics; that the factor $1-y$ ``cuts off'' the critical point which contributed to the dominant asymptotics in the negative drift case gives an analytic reason for why the combinatorial factor of drift affects the asymptotic growth rate.  It remains only to find the points in $\mV$ with the same coordinate-wise modulus and determine the contribution of each.  A quick calculation shows that there are at most four such points, the subset of $(x,y,t)=(x,1,\pm S(x,1))$ such that $x=\pm1$ and $|t|=|\mS|$. Calculating the contributions of each gives the asymptotics of $C(1,1,t)$ which appear in Table~\ref{tab:guessed}.

\subsection{The Sporadic Cases}
\[ \diag{N,W,SE} \qquad \diag{E,SE,W,NW} \qquad \diag{NW,SE,N,S,E,W} \]
Aside from the four algebraic models which were discussed in the opening section, for which asymptotics are already known, there are only three models in Table~\ref{tab:guessed} whose asymptotics remain to be proven.  The asymptotics are proven case by case, and each is a straightforward application of analytic combinatorics in several variables.  The asymptotics of $C(1,1,t)$ for these models were originally proven by Bousquet-M{\'e}lou and Mishna~\cite{BoMi10} via different methods (our results on the asymptotics of boundary returns for the first and third sporadic models, in the next section, are not covered by their work).

\begin{example}
\label{ex:E1}
For the first step set $\mS = \{N,W,SE\}$, Equation~\eqref{eq:diag} simplifies to 
\[ C(1,1,t) = \Delta \left( \frac{(x^2-y)(1-\ox\oy)(x-y^2)}{(1-x)(1-y)(1-xyt(\oy+y\ox+x))} \right).\]
As $x^2-y = (x-1)(x+1)-(y-1)$ one can re-write the rational function above as the sum of two rational functions with denominators $(1-x)(1-xyt(\oy+y\ox+x))$ and $(1-y)(1-xyt(\oy+y\ox+x))$, simplifying the singular geometry to aid calculations.  Each of the summands admits three minimal critical points:
$\bp_1 = (1,1,1/3), \bp_2 = (\nu,\nu^2,\nu^2/3),$ and $\bp_3 = (\nu^2,\nu,\nu/3)$, 
where $\nu = e^{2\pi i/3}$.  Only $\bp_1$ turns out to affect dominant asymptotics: it gives a contribution of $3^n\cdot n^{-3/2}\cdot \frac{3\sqrt{3}}{4\sqrt{\pi}}+ O(3^n\cdot n^{-5/2})$
at each of the two summands, while $\bp_2$ and $\bp_3$ give contributions of $O(3^n/n^2)$.
\end{example}
\vspace{0.02in}

\begin{example}
For the second step set $\mS = \{E,SE,W,NW\}$, Equation~\eqref{eq:diag} simplifies to
\[ F(t) = \Delta \left( \frac{(x+1)(\ox^2-\oy)(x-y)(x+y)}{1-xyt(x+x\oy+y\ox+\ox)} \right). \]
This case turns out to be easy to analyze, since the denominator is smooth.  There are two points which satisfy the critical point equations: $\bp_1 = (1,1,1/4)$ and $\bp_2 = (-1,1,1/4)$, both of which are minimal and smooth.  As the numerator has a zero of order 2 at $\bp_1$ but order 3 at $\bp_2$, in fact only $\bp_1$ contributes to the dominant asymptotics. 
\end{example}
\vspace{0.02in}

\begin{example}
For the final step set $\mS = \{NW,SE,N,S,E,W\}$, Equation~\eqref{eq:diag} simplifies to
\[ C(1,1,t) = \Delta \left( \frac{(x-y^2)(1-\ox\oy)(x^2-y)}{(1-x)(1-y)(1-txy(x+y+x\oy+y\ox+\ox+\oy))} \right). \]
Now $\bp = (1,1,1/6)$ is the only critical minimal point, and the analysis at $\bp$ is analogous to Example~\ref{ex:E1}.
\end{example}

\section{Walks Returning to the Boundary}
\label{sec:Boundary}

In the previous section we derived asymptotics for $[t^n]C(1,1,t)$ -- that is, the number of walks ending anywhere.  Combinatorially, it is also of interest to determine asymptotics for the number of walks returning to the origin or one of the bounding coordinate axes.  Theorem~\ref{thm:diag} gives a simple link between the rational functions whose diagonals we take to get these sequences: one simply multiplies the rational function $F(x,y,t)$ whose diagonal determined $C(1,1,t)$ by $1-x$, $1-y$, or both factors to get the appropriate representation.  This allows one to see immediately that for the highly symmetric and negative drift models the contributing points determining dominant asymptotics are unchanged, meaning the exponential growth of all walks and walks returning to either/both axes are the same (the order of vanishing of the numerator at these points will increase, however, meaning the polynomial growth term will be changed, along with the constant).  For the positive drift models, in contrast, the point $(x_1,y_1,t_1)$ (along with the other points in $\mV$ with the same modulus) will determine dominant asymptotics for walks returning to the $x$-axis or origin as the factor of $1-y$ in the numerator will cancel with the one present in the denominator.  This means that the exponential growth of walks returning to the $x$-axis or origin for positive drift models will be less than the exponential growth for the total number of walks, which itself is equal to the exponential growth for the number of walks returning to the $y$-axis.  The complete list of these asymptotics is given in Tables~\ref{tab:boundary1} and~\ref{tab:boundary2}, and we note that for each of the algebraic models 5--8 the minimal polynomial of $C(x,y,t)$ has been determined previously~\cite{ BoMi10, BoKa10}, meaning the asymptotics corresponding to these 4 models are already known.  These results prove numerically guessed asymptotics of Bostan et al.~\cite{BCHKP15+}.

\begin{table}[ht]
{\footnotesize
\centering
\begin{tabular}{ | c | c | c | c @{ \hspace{0.01in} }@{\vrule width 1.2pt }@{ \hspace{0.01in} }  c | c | c | c |  }
  \hline
   $\mS$ & $C(0,1,t)$ & $C(1,0,t)$ & $C(0,0,t)$ & $\mS$ & $C(0,1,t)$ & $C(1,0,t)$ & $C(0,0,t)$ \\ \hline
  &&&&&&& \\[-5pt] 
  \diag{N,S,E,W}  & $\frac{8}{\pi} \cdot \frac{4^n}{n^2}$ & 
  $\frac{8}{\pi} \cdot \frac{4^n}{n^2}$ &  
  $\delta_n \frac{32}{\pi} \cdot \frac{4^n}{n^3}$ &
  \diag{NE,SE,NW,SW}  & $\delta_n \frac{4}{\pi} \cdot \frac{4^n}{n^2}$ & 
  $\delta_n \frac{4}{\pi} \cdot \frac{4^n}{n^2}$ &  
  $\delta_n \frac{8}{\pi} \cdot \frac{4^n}{n^3}$ \\
  \diag{N,S,NE,SE,NW,SW} &$\frac{3\sqrt{6}}{2\pi} \cdot \frac{6^n}{n^2}$ & 
  $\delta_n \frac{2\sqrt{6}}{\pi} \cdot \frac{6^n}{n^2}$ &  
  $\delta_n \frac{3\sqrt{6}}{\pi} \cdot \frac{6^n}{n^3}$ &
  \diag{N,S,E,W,NW,SW,SE,NE}  & $\frac{32}{9\pi} \cdot \frac{8^n}{n^2}$ & 
  $\frac{32}{9\pi} \cdot \frac{8^n}{n^2}$ &  
  $\frac{128}{27\pi} \cdot \frac{8^n}{n^3}$  \\
  \diag{NE,NW,S}  & $\frac{3\sqrt{3}}{4\sqrt{\pi}}  \frac{3^n}{n^{3/2}}$ & 
  $\delta_n \frac{4\sqrt{2}}{\pi}  \frac{(2\sqrt{2})^n}{n^2}$ & 
  $\epsilon_n \frac{16\sqrt{2}}{\pi} \frac{(2\sqrt{2})^n}{n^3}$ &
  \diag{N,NW,NE,S}  & $\frac{8}{3\sqrt{\pi}} \frac{4^n}{n^{3/2}}$ & 
  $\delta_n \frac{4\sqrt{3}}{\pi}  \frac{(2\sqrt{3})^n}{n^2}$ & 
  $\delta_n \frac{12\sqrt{3}}{\pi}  \frac{(2\sqrt{3})^n}{n^3}$ \\
  % %%%%%%%%%%%%%
  \diag{NE,NW,E,W,S} & $\frac{5\sqrt{10}}{16\sqrt{\pi}}  \frac{5^n}{n^{3/2}}$  & 
  $\frac{\sqrt{2}A^{3/2}}{\pi}  \frac{(2A)^n}{n^2}$ & 
  $\frac{2A^{3/2}}{\pi}  \frac{(2A)^n}{n^3}$ &
  \diag{N,NE,NW,SE,SW}  & $\frac{5\sqrt{10}}{24\sqrt{\pi}}  \frac{5^n}{n^{3/2}}$ & 
  $\delta_n \frac{4\sqrt{30}}{5\pi}  \frac{(2\sqrt{6})^n}{n^2}$ & 
  $\delta_n \frac{24\sqrt{30}}{25\pi}  \frac{(2\sqrt{6})^n}{n^3}$ \\
  % %%%%%%%%%%%%%
  \diag{N,NW,NE,E,W,S} & $\frac{\sqrt{3}}{\sqrt{\pi}} \frac{6^n}{n^{3/2}}$ &
  $\frac{2\sqrt{3}B^{3/2}}{3\pi} \frac{(2B)^n}{n^2}$ & 
  $\frac{2B^{3/2}}{\pi} \frac{(2B)^n}{n^3}$ &
  \diag{N,E,W,NE,NW,SE,SW} & 
  $\frac{7\sqrt{21}}{54\sqrt{\pi}} \frac{7^n}{n^{3/2}}$ & 
  $ \frac{D}{285\pi} \frac{(2K)^n}{n^2}$ & 
  $\frac{2E}{1805\pi} \frac{(2K)^n}{n^3}$ \\ 
  % %%%%%%%%%%%%%%%
  \diag{N,W,SE}   & $\frac{27\sqrt{3}}{8\sqrt{\pi}} \cdot \frac{3^n}{n^{5/2}}$ & 
  $\frac{27\sqrt{3}}{8\sqrt{\pi}} \cdot \frac{3^n}{n^{5/2}}$ & 
  $\sigma_n \frac{81\sqrt{3}}{\pi} \cdot \frac{3^n}{n^4}$ &
  \diag{E,SE,W,NW}  & $\delta_n \frac{32}{\pi} \cdot \frac{4^n}{n^3}$ & 
  $\frac{32}{\pi} \cdot \frac{4^n}{n^3}$ & 
  $\delta_n \frac{768}{\pi} \cdot \frac{4^n}{n^5}$ \\
  % %%%%%%%%%%%%%
  \diag{NW,SE,N,S,E,W} & $\frac{27\sqrt{3}}{8\sqrt{\pi}} \cdot \frac{6^n}{n^{5/2}}$ & 
  $\frac{27\sqrt{3}}{8\sqrt{\pi}} \cdot \frac{6^n}{n^{5/2}}$ & 
  $\frac{27\sqrt{3}}{\pi} \cdot \frac{6^n}{n^4}$ & &&&\\
\hline
\end{tabular}
\vspace{0.05in}

\caption{Asymptotics of boundary returns for the highly symmetric, positive drift, and sporadic cases. } \label{tab:boundary1}

\begin{tabular}{ | c | c @{ \hspace{0.01in} }@{\vrule width 1.2pt }@{ \hspace{0.01in} } c | c |  }
  \hline
   $\mS$ & $C(0,1,t)$ & $\mS$ & $C(0,1,t)$ \\ \hline
  &&& \\[-5pt] 
  \diag{N,SE,SW}  & $
    \left(\epsilon_n \frac{448\sqrt{2}}{9\pi}  +   
    \epsilon_{n-1} \frac{640}{9\pi}  +
    \epsilon_{n-2} \frac{416\sqrt{2}}{9\pi} + 
    \epsilon_{n-3} \frac{512}{9\pi}  
    \right)\cdot \frac{(2\sqrt{2})^n}{n^3} 
   $ & 
  \diag{S, SW, SE, N}  & $\left(\delta_n \frac{36\sqrt{3}}{\pi} + \delta_{n-1} \frac{54}{\pi} \cdot \frac{(2\sqrt{3})^n}{n^3}\right)\cdot \frac{(2\sqrt{3})^n}{n^3}$ \\
  % %%%%%%%%%%%%%
  \diag{SE, SW, E, W, N} & $\frac{4A^{7/2}}{\pi} \cdot \frac{(2A)^n}{n^3}$  &
  \diag{S, SE, SW, NE, NW}  & $\left( \delta_n \frac{72\sqrt{30}}{5\pi} + \delta_{n-1}\frac{864\sqrt{5}}{25\pi} \right) \cdot \frac{(2\sqrt{6})^n}{n^3}$\\
  % %%%%%%%%%%%%%  
  \diag{S, SW, SE, E, W, N} & $\frac{3B^{7/2}}{2\pi} \cdot \frac{(2B)^n}{n^3}$  &
  \diag{NE, NW, E, W, SE, SW, S} & 
  $\frac{6(4571+1856\sqrt{6})\sqrt{23-3\sqrt{6}}}{1805\pi} \cdot \frac{(2K)^n}{n^3}$ \\[+1mm]
  % %%%%%%%%%%%%%%%
\hline
\end{tabular}
\vspace{0.05in}

\caption{Asymptotics of $C(0,1,t)$ for negative drift cases; other asymptotics of $\mS$ are the same as those of $-\mS$ above.} \label{tab:boundary2}
\vspace{-0.2in}

\[ A = 1+\sqrt{2}, \qquad B = 1+\sqrt{3}, \qquad K = 1+\sqrt{6}, \qquad D = (156+41\sqrt{6})\sqrt{23-3\sqrt{6}}, \qquad E = (583+138\sqrt{6})\sqrt{23-3\sqrt{6}}\]
\begin{center}
$\delta_n = 1$ if $n \equiv 0$ mod 2, $\sigma_n = 1$ if $n\equiv0$ mod 3, and $\epsilon_n=1$ if $n \equiv 0$ mod $4$ -- each is 0 otherwise
\end{center}
}
\end{table}

%---------------------------------------
\section{Conclusion}
\label{sec:Conc}
%---------------------------------------

Lattice paths restricted to the non-negative quarter plane are well-studied objects, and there are many approaches to their enumeration.  In this article we have tried to highlight the benefits of using diagonals of rational functions: compact representations of GFs, effective determinations of asymptotics, and clear links between analytic and combinatorial properties.  After determining asymptotics of the 23 D-finite models in the quarter plane there are clear generalizations still left to be worked on: models with longer or weighted steps, models in higher dimensions, and restrictions to different regions or other lattices are a few obvious ones.  Current work is ongoing in several of these areas, however in some situations the diagonal representations obtained by the kernel method do not have a singular variety which admits critical minimal points, meaning a deeper singularity analysis is required.

%--------------------------------------
% bilbiography
%--------------------------------------
%\footnotesize
%\bibliographystyle{plain}
%\bibliography{bibl}

\let\oldbibliography\thebibliography
\renewcommand{\thebibliography}[1]{\oldbibliography{#1}
\setlength{\itemsep}{-2pt}} %Reducing spacing in the bibliography.

\end{document}